\documentclass[11pt]{article}
\usepackage{times}
\usepackage{graphicx}
\usepackage{amsmath,amsthm,amsfonts}
\usepackage{epstopdf}
\usepackage[height=8.25in,width=5.5in]{geometry}
\theoremstyle{definition}
\newtheorem{example}{Example}

\newcommand{\N}{\mathbb{N}} 
\newcommand{\Q}{\mathbb{Q}}

\newcommand{\alg}{\mathcal{A}}
\newcommand{\ra}{\rightarrow}

\title{Counterexamples in the Theory of Fair Division}
\author{Theodore P. Hill \footnote{School of Mathematics, Georgia Institute of Technology, Atlanta, GA 30332. \newline \tt{hill@math.gatech.edu}} \ 
\footnote{Department of Mathematics, California Polytechnic State University, San Luis Obispo, CA 93407. \newline \tt{kmorriso@calpoly.edu}}
 \and Kent E. Morrison \footnotemark[\value{footnote}]
}
\date{December 22, 2008}
\begin{document}
\maketitle

\large
\renewcommand{\baselinestretch}{1.25}   
\normalsize

The formal mathematical theory of fair division has a rich history dating back at least to Steinhaus in the 1940's \cite{Ste1,Ste2,Ste3,Ste4}. In recent work in this area, several general classes of errors have appeared along with confusion about the necessity and sufficiency of certain hypotheses. It is the purpose of this article to correct the scientific record and to point out with concrete examples some of the pitfalls that have led to these mistakes. These examples may serve as guideposts for future work. 

As a basic starting point for discussing fair division, consider the classical Moving-Knife procedure below, a continuous version of a discrete procedure first described by Steinhaus \cite{Ste3} and attributed to Banach and Knaster:

\textbf{Moving-Knife Procedure} \cite[p.2]{DSp}. ``A knife is slowly moved at constant speed parallel to itself over the top of the cake. At each instant the knife is poised so that it could cut a unique slice of the cake. At time goes by the potential slice increases monotonely from nothing until it becomes the entire cake. The first person to indicate satisfaction with the slice then determined by the position of the knife receives that slice and is eliminated from further distribution of the cake. (If two or more participants simultaneously indicate satisfaction with the slice, it is given to any of them.) The process is repeated with the other $n-1$ participants and with what remains of the cake."

The only implicit assumptions here are that the values are nonnegative, additive, and, in the direction the knife is moving, continuous. The knife need not even be perfectly straight and need not be moved perfectly parallel to itself, and the cake need not be simply connected (as traditional angel-food cakes, with a hole through the center, are not), nor even connected (the cake could have been dropped on the table and landed in pieces). Even though written informally, the ideas and conclusions are clear and correct. 

In some subsequent articles about fair division, however, even with the appearance of greater mathematical rigor, these implicit hypotheses have been altered, sometimes leading to faulty conclusions.

\textbf{FRAMEWORK}

The basic setting in fair-division theory is that there is an object (often called the \emph{cake}) that is to be divided among $n$ players who may have different relative values of the cake. In general the cake, here denoted by $X$, may be an arbitrary non-empty set, but often $X$ is assumed to be a three-dimensional object (e.g. \cite{Str}), a rectangle (e.g. \cite{BT}), or even the unit interval \cite{BJK}.

The \emph{allowable portions} of the cake consist of a non-empty collection of subsets of $X$, here denoted $\mathcal{A}$, which is often assumed to be either an algebra (containing the empty portion $\emptyset$ , and closed under finite unions and complements) or a sigma algebra (algebra closed under countable unions). In some cases, such as when $X$ is a finite set or the values are finitely additive, as in Example 1 below, $\mathcal{A}$ may consist of all subsets of $X$. In other cases, $\mathcal{A}$ may be a proper set of subsets of $X$, such as ``vertical" slices of a rectangular cake $X$, or finite unions of open intervals of a one-dimensional cake $X$, or the Borel subsets of $X$ (smallest sigma algebra containing the open intervals or balls of an $n$-dimensional $X$).

A player's \emph{value} (or value measure) $v$ is a function $v:\mathcal{A} \ra [0,1]$ satisfying $v(\emptyset)=0$ and $v(X)=1$; the number $v(A)$ represents the relative value of portion $A \in \mathcal{A}$ of the cake $X$ to a player with value $v$. Cake-cutting problems usually assume that the cake is either infinitely divisible ($v(A) \in (0,1) \Rightarrow  \mbox{ there exists } B \subset A,  B \in \mathcal{A} \mbox{ with }\, 0 < v(B) < v(A)$, or defined via an equivalent bisection property) (e.g. [DSp]), or is purely discrete ($v(A)=\sum_{x \in A} v(x)$ when $X$ is finite or countable) (e.g. [DH]).

\textbf{ADDITIVITY}

The players' values are usually either assumed to be finitely additive ($v(A \cup B)=v(A)+v(B)$ if $A \cap B = \emptyset$) or countably additive ($v \left( \bigcup_{n=1}^\infty A_n \right)=\sum_{n=1}^\infty v(A_n)$ if $\{A_n\}$ are mutually disjoint), although some progress has also been made on the sub-additive case (e.g. \cite{DM}). For many fair division procedures such as classical cut-and-choose or the moving-knife method above, finite additivity suffices. If the players' values are further assumed to be countably additive, then it is possible to ``cut the cake in an infinite sequence of cuts without losing value" \cite[p. 20]{RW}, but if only finite additivity is assumed, some infinite procedures which have been proposed may fail.

\begin{example} 
$X$ is the unit interval $[0,1]$, and $\alg$ is the set of all subsets of $X$. Let $v$ be a finitely additive probability measure on $(X,\alg )$ that is concentrated on the rationals in $X$
(i.e., $v(\Q \cap X)=1$) satisfying $v([a,b])=b-a$ for all $0 \leq a \leq b \leq 1$. (It follows from the Hahn-Banach theorem that many such finitely additive probability measures exist, c.f \cite{DSa}.) Letting $(r_n)$ be an enumeration of $\Q \cap X$, note that $v(r_n)=0$ for all $n \in \N$, so $v$ is not countably additive since $v(\Q \cap X)=1 > 0 = \sum_{n \in \N} v(r_n)$ . To see that mass may disappear in an infinite sequence of cuts, let $A_n=[0,1/2+1/n] \setminus \{r_1,\ldots,r_n\}$ for all $n \in \N$. Since $v(r_n)=0$ for all $n$, $v(A_n)=1/2 +1/n$ for all $n$, so the $v$-values of the decreasing portions $A_n$ are strictly decreasing to $1/2$. But after an infinite number of steps, $\bigcap_{n=1}^\infty A_n = [0,1/2] \setminus \Q$, so $v\left( \bigcap_{n=1}^\infty A_n \right) \leq v(X \setminus \Q) =0$, and the limiting value of decreasing portions that have value greater than $1/2$ at every step is now zero. This shows that in assuming only finite additivity \cite[p. 10]{BT}, the infinite ``trimming procedure" in \cite[Ch. 7.3]{BT} and \cite{COMAP} does not always produce an evy-free or even fair allocation, contrary to the claim in \cite[p 38]{BT}.
\end{example}

\textbf{CONTINUITY}

In addition to hypotheses about the additivity of the underlying value measures, different basic assumptions concerning the continuity of the measures can be critical in cake-cutting. Asssuming that countably additive  measures are \emph{atomless} (or non-atomic) (e.g. \cite{Jml}) means that no single points have positive value, but even this continuity assumption requires care. If the cake is the closed unit square, for example, and the players' value measures are all uniform on the top frosting (i.e., on $\{(x,1) : 0 \leq x \leq 1 \}$), then those measures are all atomless, but the moving-knife procedure applied with the knife horizontal and moving from top to bottom does not provide a fair division since every player says ``stop" immediately. (Moving the knife at a different angle, however, does guarantee an allocation where every player gets a portion he values exactly $1/n$, and for general atomless measures, there is always a direction where the moving-knife algorithm generates a fair division, and in fact, if the direction is chosen uniformly at random, this will happen with probability 1 \cite[Corollary 2]{Jml}. The hypothesis of the measures being atomless is key to the application of many key tools in cake-cutting, especially Lyapounov's Convexity Theorem, and without non-atomicity, fair allocations do not always exist.

A stronger continuity hypothesis than non-atomicity is that the value measures are \emph{absolutely continuous}, that is, absolutely continuous with respect to Lebesgue measure (e.g. \cite{Jma}), which simply says that every portion of an $n$-dimensional cake that has no volume (i.e., its $n$-dimensional Lebesgue measure is zero) has no value to any of the players. For example, if the measures are absolutely continuous, the moving-knife method works at every angle.

 Some cake-cutting results require that the measures are not only absolutely continuous, but also mutually absolutely continuous, or mutually absolutely continuous with respect to Lebesgue measure. This may seem like an innocent technicality of nonnegativity-versus-positivity, but in fair division such differences often imply important philosophical changes in the problem. For example, if the cake is an inhomogeneous mixture of various ingredients including chocolate, then absolute continuity may be a reasonable hypothesis, but it is not one that is standard in classical or modern fair division theory (e.g., \cite{B, DSp, EHK, Jml}). But further requiring the values to be mutually absolutely continuous means that \emph{if one player likes chocolate, all the other players must like chocolate as well.} And requiring that the measures be mutually absolutely continuous with respect to Lebesgue measure (i.e., the measures are determined by probability density functions that are strictly positive almost everywhere), means that \emph{every player must like chocolate} (and must like every other part of the cake). 

In some articles, the basic standard terminology concerning continuity of measures is misleading, confusing or incorrect. In \cite[p. 1315]{BJK} for example, the authors ``postulate that the players have continuous value functions...and their measures are finitely additive." The authors then state ``we assume that the measures of the players are absolutely continuous, so no portion of cake is of positive measure for one player and zero measure for another player" (see also \cite[p. 291]{Br}). Absolute continuity, in its common usage, implies \emph{countable} additivity, and absolute continuity does not imply mutual absolute continuity, nor vice versa, so it is not clear what additivity hypotheses the authors intended.  Moreover, the main results in \cite{BJK} and \cite[Chapter 13]{Br} are not true without further assuming mutual absolute continuity with respect to Lebesgue measure (see Examples 2, 3, 4, 5 below and \cite{M}).

Even if the cake is the unit interval, continuous countably additive value measures may not have unique medians, and this may invalidate some procedures. For example, the  cake-cutting procedure Surplus Procedure (SP) is not well defined without an additional assumption guaranteeing unique medians \cite{M}, since in SP,  (1) player 1 and player 2 independently report these value functions to a referee, and (2) the ``referee determines the 50-50 points [medians]" (\cite[p. 292]{Br}, \cite[p. 1315]{BJK}). But if either player has a value function with multiple medians, then step (2) is not well defined.  

Without additional continuity assumptions, other procedures may also fail. For example, part (2) of the definition of the Equitability Procedure (EP) (\cite[p. 299]{Br}, \cite[p. 1315]{BJK}) assumes the existence of ``cut-points that equalize the common value that all players receive for each of the $n!$ possible assignments of pieces to the players from left to right." As the next example shows, even with absolutely continuous measures, without an additional assumption of mutual absolute continuity with respect to Lebesgue measure, such cut-points may not exist, so EP is not well defined \cite{M}.

\begin{example} 
 The cake is the unit interval. Player 1 values the cake uniformly, player 2's value is uniform on $(0,1/3)$ (i.e., his pdf is a.s. constant 3 on $(0,1/3)$ and zero elsewhere), and player 3's value is uniform on $(2/3,1)$. Then for the ordering 1-3-2 (from left to right), there do not exist two cut-points that equalize the values. If the second cut-point is in $(0, 2/3]$ then player 3 receives 0 but player 1 receives a positive amount. If the second (and hence both cut-points) are at 0, players 1 and 3 receive 0, and player 2 receives 1. If the second cut-point is in $(2/3, 1)$, then player 3 receives a positive amount, but player 2 receives 0. If the second cut-point is at 1, and the first cut-point is 0, then players 1 and 2 get 0, but player 3 gets 1. If the second cut-point is at 1 and the first is in $(0,1)$, then player 2 gets 0 but player 1 gets a positive amount. Finally, if both the first and the second cut-points are at 1, then player 1 gets 1, and both players 2 and 3 get zero. Thus it is not possible to equalize the common value that all players receive for the order 1-3-2.
\end{example}

\textbf{PARETO OPTIMALITY}
 
The Pareto optimality (also called efficiency) of an allocation of allowable portions to players concerns comparison of the given allocation with other possible allocations. If no ``better" allocations exist, the allocation is called \emph{Pareto optimal} (efficient). The key notion here, however, is the interpretation of ``better", and there are two different meanings meanings that are in common use. For example, in \cite{Jma} the author defines efficiency as there being no other allocation ``that gives both players more of the cake" \cite[p. 278]{Jma}, a ``weak" Pareto optimality. In \cite[p. 289]{Br} and \cite[p. 1314]{BJK}, however, efficiency is defined as ``There is no other allocation that is better for one person and at least as good for the other" \cite[p. 1314]{BJK}, a ``strong" Pareto optimality. Without additional continuity assumptions, the main result in \cite{Jma} holds for weak Pareto optimality (as is correctly stated), but fails for strong Pareto optimality. In \cite{BJK} the authors state that in problems of fair division of a divisible good, ``the well-known 2-person, 1-cut cake-cutting procedure `I cut, you choose'" is (strong) Pareto-optimal. Even assuming absolute continuity of the value measures, however, cut-and-choose is not even weak Pareto optimal.
 
\begin{example} 
The ``cake" is the unit square, and player 1 values only the top half of the cake and player 2 only the bottom half (and on those portions, the values are uniformly distributed). If player 1 is the cutter and cuts vertically, his uniquely optimal cut-and-choose solution is to bisect the cake exactly, in which case each player receives a portion he values exactly 1/2. But an allocation of the top half of the cake to player 1, and the bottom half to player 2 gives both players a portion each values 100\%, which is strictly better for both players. Thus cut-and-choose is neither strong nor weak Pareto optimal in general. If player 1 cuts horizontally, on the other hand, his uniquely optimal risk-adverse cut-and-choose point is  the line  y = 3/4, in which case he receives a portion he values at 50\% of the cake, and player 2 chooses the bottom portion and receives a portion he values at 100\% of the cake.  But an allocation of the top half of the cake to player 1, and the bottom half to player 2 is at least as good for player 2 in both cases, and is strictly better for player 1. Thus if player 1 cuts horizontally, cut-and-choose is weak Pareto optimal but not strong Pareto optimal.
\end{example}

In \cite[p. 299]{Br} and \cite[p. 1318]{BJK}, the authors claim that EP is (strong) Pareto optimal, and in \cite[p. 1320]{BJK}, that SP is Pareto optimal. Both those claims are incorrect, even when all the value measures are represented by probability density functions that are strictly positive everywhere. The underlying reason is that both EP and SP allocate contiguous portions to each player, and as noted in \cite[p. 149]{BT}, ``satisfying contiguity may be inconsistent with satisfying efficiency." This is illustrated in the next two examples, which show that EP and SP, respectively, are not in general Pareto optimal.

\begin{example} 
The cake $X$ is the unit interval $[0,1]$, and there are two players. Player 1's value function is 1.6 on $ (0, 1/4)$ and on $(1/2, 3/4)$ and is 0.4 elsewhere; and player 2's is 1.6 on $(1/4, 1/2)$ and $(3/4, 1)$ and 0.4 elsewhere. Then SP cuts the cake at $1/2$, and each player receives a portion worth exactly 0.5. But allocating $(0, 1/4)$ and $(1/2, 3/4)$ to player 1, and the rest to player 2, gives each player a portion he values 0.8, which is strictly better for each player than the SP allocation. Thus SP is not even weak Pareto optimal.
\end{example}

\begin{example} 
The cake is the unit interval $[0,1]$, and there are three players. Player 1's value function is 2.4 on $(0, 1/6)$ and on $(1/2, 2/3)$ and is 0.3 elsewhere; player 2's value function is 2.4 on $(1/6, 1/3)$ and $(2/3, 5/6)$ and 0.3 elsewhere; and player 3's is 2.4 on $(1/3, 1/2)$ and $(5/6, 1)$ and 0.3 elsewhere. Then one EP solution allocates $(0, 1/3)$ to player 1, $(1/3, 2/3)$ to player 3, and the rest to player 2, and each player receives a portion worth exactly $0.45$. But allocating $(0, 1/6)$ and $(1/2, 2/3)$ to player 1, $(1/6, 1/3)$ and $(2/3, 5/6)$ to player 2, and the rest to player 3 gives each player a portion he values 0.8, which is strictly better for each player than the EP allocation, so EP is not Pareto optimal. 
\end{example}

\textbf{LIMITED CUTS}

EP for 3 players makes an allocation based on two cuts, whereas if five cuts are allowed in Example 5, each player may receive a portion he considers 80\% of the value. Whether allocations should be limited to $n-1$ cuts in cases like this where a few more cuts will permit a much better solution, is debatable. In some procedures such as cut-and-choose, the allocation is effected by $n-1$ cuts in a fair-division problem with $n$ players. In \cite[p. 297]{Br} and \cite[p. 1318]{BJK} the authors state ``an envy-free allocation that uses $n-1$ parallel, vertical cuts is always efficient" (i.e., strong Pareto optimal) without additional hypotheses. However, both those claims and the corresponding Proposition 7.1 of \cite[pp. 149--150]{BT}, which assumed the setting of \cite{BDT} (atomless or absolutely continuous measures), are incorrect. The next example illustrates this using the classical envy-free procedure cut-and-choose.

\begin{example} 
The cake $X$ is the unit interval; player 1 values it uniformly, and player 2 values only the left- and right-most quarters of the interval, and values them equally and uniformly (i.e., the probability density function representing player 1's value is constant 1 on $[0,1]$, and that of player 2 is constant 2 on $[0, 1/4]$ and on $[3/4,1]$, and zero otherwise.) If player 1 is the cutter, his unique cut-point is at $x = 1/2$, and each player will receive a portion he values at exactly $1/2$. The allocation of the interval $[0, 1/4]$ to player 2 and the rest to player 1, however, gives player 1 a portion he values $3/4$, and player 2 a portion he values $1/2$ again, so cut-and-choose (which is an envy-free allocation for two players) is not weak Pareto optimal, not even among 1-cut allocations. 
\end{example}

Similarly, Corollary 7.1 of \cite[p. 151]{BT}, which says that Stomquist's well known envy-free moving-knife procedure for three players [Str]  is ``C-efficient" (strong Pareto optimal among allocations using $n-1$ cuts) , is also incorrect without additional assumptions, as is shown in the next example.

\begin{example} 
 $X$ is the closed unit interval $[0,1]$, and there are three players. Player 1 values $X$ uniformly on $X$, and player 2 and player 3 both value $X$ uniformly on $(.4,.6)$. Then the Stromquist 4-knife envy-free procedure, moving from left to right, gives portion (0,1/3) to player 1, and gives players 2 and 3 portions $(1/3, 1/2)$ and $(1/2,1)$; it does not matter who gets which. Thus the Stromquist procedure allocates a portion to player 1 which he values $1/3$, while players 2 and 3 each gets a portion he values $1/2$. But the allocation $(0,.4)$ to player 1, $(.4,.5)$ to player 2 and the rest to player 3 is strictly better for player 3 and equally good for the other two, so the $n -1$ cut Stromquist procedure is not strong Pareto optimal. 
\end{example}

\textbf{INCENTIVE COMPATIBILITY}

The complex and somewhat game-theoretic issues of what information is available to the players, whether players must act truthfully, and whether the referee (if the procedure requires one) must act in a prescribed manner (e.g. if several fair solutions are available), is crucial and has also been a source of errors. An allocation procedure is sometimes called \emph{incentive compatible} if all of the (risk-adverse) players fare best when they truthfully reveal any private information asked for by the procedure. Under essentially every cake-cutting allocation procedure, a player who does not act truthfully risks receiving a share he values less than $1/n$ of the cake, as is illustrated in the next example. Thus the second parts of the conclusions in \cite[Proposition 13.3]{Br} and \cite[Theorem 3]{BJK}, namely, that under procedure EP a player who is not truthful may receive less than $1/n$ of the cake, and the statement ``maximin strategies under SP require each player to be truthful" \cite[p. 290]{Br}, are vacuous.

\begin{example}  
$X$ is the closed unit interval $[0,1]$, there are two players, and both report (falsely) that their value measures are uniform on $[0,1]$. Based on these reported values, each allocation procedure assigns some subset $A$ of $X$ to player 1, and the complement $X \setminus A$ to player 2 (e.g., $A = [0,1/2]$). But if player 2's true value of $A$ and player 1's true value of $X \setminus A$ are both 1, then both players have received a portion valued at zero. Thus reporting their measures falsely resulted in both receiving a portion less than $1/n$.
\end{example}

There are a variety of notions of incentive compatibility.  In \cite[p. 294]{Br} and \cite[p. 1316]{BJK}, for example, the authors define an allocation procedure to be \emph{strategy-vulnerable} if a ``player can, by misrepresenting its value function, assuredly do better, whatever the value function of the other player"; and otherwise the procedure is called \emph{strategy-proof}. But under this definition, as the next easy example shows, \emph{every} fair procedure is strategy-proof, so without additional assumptions the second conclusion of \cite[Proposition 13.1]{Br} and
\cite[Theorem 1]{BJK} that ``any procedure that makes $e$ the cut-point is strategy vulnerable" is false in general \cite{M}.

\begin{example}  
 The cake is the closed unit interval, and there are $n$ players all of whom have identical values (e.g. all value it uniformly). Since every fair procedure allocates disjoint subsets of the cake, one to each player, at least one of the players will receive a portion he values no more than $1/n$. That player will not do ``assuredly better" than his fair share of $1/n$. Hence every fair allocation procedure is strategy-proof, including the one that makes $e$ the cut-point, contradicting the conclusion of Theorem 1. This example also shows that the first conclusion of \cite[Theorem 1]{BJK} (that SP is strategy-proof)  and \cite[Theorem 2]{BJK} (that EP is strategy-proof) are trivial, since both SP and EP are fair.
\end{example}

In \cite{M}, Brams, Jones and Klamler acknowledge the error in \cite[Theorem 1]{BJK} and \cite[Proposition 13.1]{Br} concerning strategy-proof, but they attempt to correct it by saying that it can be rectified by substituting ``do at least as well and sometimes better" for ``do assuredly better" in their definition of strategy-vulnerability. That claim is also incorrect, as the next example, a modification of an example of K. Teh, shows.

\begin{example} 
 The cake $X$ is the unit interval, there are two players, and both players' value measures are represented by probability densities that are positive almost everywhere. Let player 1's true measure be given by the cumulative distribution function $F$ (so $F(x)$ is the value player 1 places on the interval $[0,x]$, and $F$ is thus strictly increasing everywhere, since $dF/dx > 0$ everywhere). This implies that $F$ has a unique median at $a$. Suppose that rather than reporting his true cdf $F$, player 1 reports the cdf $F^*$ given by
 \[  F^*(x) =
           \begin{cases}
              \frac{1}{2} -2 (\frac{1}{2} - F(x))^2   &\text{if $x \leq a$} \\
              \frac{1}{2} +2 (\frac{1}{2} - F(x))^2   &\text{if $x \geq a$}
           \end{cases}
\]
Then if player 2 reports the same median $a$ as player 1, the SP cut-point is at $a$, and both players receive exactly half the cake. If player 2 reports a median $b >  a$, then player 1 will receive a strictly larger portion by reporting $F^*$ instead of his true $F$ if and only if the cut-point $c$ in $(a,b)$ satisfies
\[   \frac{F^*(c)-1/2}{F^*(b)-1/2}  < \frac{F(c)-1/2}{F(b)-1/2} .
\]
But this follows from the definition of $F^*$ since
\[   \frac{F^*(c)-1/2}{F^*(b)-1/2} = \frac{2(F(c)-1/2)^2}{2(F(b)-1/2)^2} 
       < \frac{F(c)-1/2}{F(b)-1/2},
\]
since $c < b$ and $F$ is strictly increasing.

The argument for $b < a$ is similar, but in this case player 1 receives the rightmost interval and the cut-point $c$ moves to the left, giving player 1 an even greater portion. Thus, no matter what value measure player 2 reports to the referee, if player 1 reports $F^*$, then player 1 will do  ``at least as well and sometimes better", so SP is not strategy-proof in the revised sense either. (Of course, if player 2 lies even more than player 1 does, then player 1 may do worse  than he would have if both players had reported their true measures.)
\end{example}

Some of the errors identified above may serve as a challenge and opportunity for students and mathematicians to develop and make rigorous. As can be seen in the moving-knife procedure above, it is often possible to express practical, yet clean, clear, and beautiful logical conclusions without using highly technical language.

 \end{document}